\documentclass[11pt]{article}

\usepackage[utf8]{inputenc}
\usepackage[T1]{fontenc}
\usepackage{lmodern}
\usepackage{microtype}
\usepackage{amsmath,amssymb,amsthm,mathtools}
\usepackage{enumitem}
\usepackage[margin=1.08in]{geometry}
\usepackage[hidelinks]{hyperref}
\allowdisplaybreaks

\newcommand{\N}{\mathbb N}

\newcommand{\eps}{\varepsilon}
\newcommand{\Lex}[1]{\mathcal L_{#1}^{\exists}}
\newcommand{\Lall}[1]{\mathcal L_{#1}^{\forall}}
\newcommand{\Lthick}[2]{\mathcal L_{#1}^{\forall,#2}}

\newcommand{\dist}{\operatorname{dist}}

\newcommand{\Cone}{\operatorname{Cone}}

\newtheorem{theorem}{Theorem}[section]
\newtheorem{lemma}[theorem]{Lemma}
\newtheorem{corollary}[theorem]{Corollary}
\newtheorem{proposition}[theorem]{Proposition}

\theoremstyle{definition}
\newtheorem{definition}[theorem]{Definition}
\newtheorem{example}[theorem]{Example}
\theoremstyle{remark}
\newtheorem{remark}[theorem]{Remark}

\title{Cotlar Properties, Hyperbolicity, and\\
Unconditional Fourier Decompositions}
\author{Tao Mei\\Texas A\&M University\\\texttt{taomei@tamu.edu}}
\date{}

\begin{document}
\maketitle

\begin{abstract}
Let a discrete group $\Gamma$ act by isometries on a metric space $X$.  We show that actions on geodesic hyperbolic spaces satisfy a branch Cotlar property. Under a finite low-orbit-frequency hypothesis, this yields uniform $L^4$-unconditionality for branch subspaces associated with sufficiently separated subsets; when the union of these branches has finite complement, the corresponding signed projections are uniformly bounded on
$L^p(\widehat\Gamma)$ for every $1<p<\infty$.  Finally, for connected vertex-transitive graphs, we prove that a scale-wise branch Cotlar property characterizes hyperbolicity of $X$.   \end{abstract}



\section*{Introduction}

For a trigonometric polynomial
\[
 f(e^{i\theta})=\sum_{k=-N}^{N}a_ke^{ik\theta},
\]
let $P_{>0}$ and $P_{<0}$ denote the projections onto the strictly positive
and strictly negative Fourier modes.  The classical Hilbert transform
\[
 H=-iP_{>0}+iP_{<0},\qquad H1=0,
\]
is bounded on $L^p(\mathbb T)$ for every $1<p<\infty$.  Mei and Ricard
constructed a noncommutative analogue on free group von Neumann algebras by
splitting a free group into its first-letter branches and proved uniform
$L^p$-bounds for all signed sums of the associated branch projections
\cite{MR17}.  A key ingredient in their argument is a Cotlar identity
adapted to the separated first-letter branches of a free group.

The natural terminology for this phenomenon is unconditionality.  A family
of closed subspaces $(X_i)_{i\in I}$ of an $L^p$-space is called
\emph{$L^p$-unconditional} if there is a constant $C_p$ such that
\begin{equation}\label{eq:intro-unconditional}
 \left\|\sum_{i\in F}\eps_i x_i\right\|_p
 \leq C_p\left\|\sum_{i\in F}x_i\right\|_p
\end{equation}
for every finite $F\subset I$, every $x_i\in X_i$, and every choice
$|\eps_i|\leq1$.  For unimodular coefficients, applying
\eqref{eq:intro-unconditional} a second time with the conjugate coefficients
gives the usual two-sided equivalence.  This is more precise here than the
classical terminology ``$\Lambda_p$-property'', which usually concerns an
$L^p$--$L^2$ estimate on one frequency set rather than unconditional signed
projections onto a family of subspaces.  Fourier multipliers and their
noncommutative analogues have attracted sustained attention; see, for
example, \cite{Pa26} and the references therein.

The full results are formulated for isometric actions on metric spaces.  To
keep the introduction concrete, we summarize the two principal theorems below
in the special case of a finitely generated group $\Gamma$ acting by left
translation on a Cayley graph
\[
 X=\operatorname{Cay}(\Gamma,S).
\]
We write $e$ for the identity, $|g|=\dist(e,g)$ for the word length, and
 $S_m(e)$ for the corresponding  sphere.  For
$A\subset X$ and $r\geq0$, let
\[
 \Lex{A}
 =\{g\in\Gamma:\text{some geodesic from $e$ to $g$ meets $A$}\},
\]
and
\[
 \Lthick{A}{r}
 =\left\{g\in\Gamma:\begin{array}{l}
 \text{there exists $a\in A$ such that every geodesic from $e$ to $g$}\\
 \text{meets the ball $B_r(a)$}
 \end{array}\right\}.
\]
Thus $\Lex{A}$ is an exact existential branch, while
$\Lthick{A}{r}$ is a universal $r$-thick branch.  Let $\lambda_g$ be the
left regular representation of $\Gamma$, and let $L^p(\widehat\Gamma)$
denote the associated noncommutative $L^p$-space.  For
$\mathcal E\subset\Gamma$, let $P_{\mathcal E}$ denote the Fourier
projection onto the closed span of $\{\lambda_g:g\in\mathcal E\}$ in
$L^2(\widehat\Gamma)$.  Given a family $\mathcal A=(A_i)_{i\in I}$ and
finitely supported coefficients $|\eps_i|\leq1$, set
\[
 \mathcal H^{\exists}_{\mathcal A,\eps}
   =\sum_i\eps_iP_{\Lex{A_i}},
 \qquad
 \mathcal H^{\forall,r}_{\mathcal A,\eps}
   =\sum_i\eps_iP_{\Lthick{A_i}{r}}.
\]
We also use the Gromov product
\[
 (g\mid h)_e=\frac{|g|+|h|-|g^{-1}h|}{2}.
\]
These are the Cayley-graph specializations of the notation introduced in
Section~1.

\paragraph{Theorem A: separated branches and unconditionality
(Theorem~\ref{thm:main} and Corollary \ref{cor:sphere}).}
Assume that $X$ is $\delta$-hyperbolic.  Let
$(A_i)_{i\in I}\subset S_m(e)$ and   $\dist(A_i,A_j)\geq 4\delta$ for $i\neq j$. 
Then the exact branch subspaces generated by $\lambda(\Lex{A_i})$ form an
$L^4$-unconditional family on their closed linear span; equivalently,
\[
 \left\|\mathcal H^{\exists}_{\mathcal A,\eps}x\right\|_4
 \lesssim\|x\|_4
\]
uniformly in the number of branches and in the coefficients, for $x$ with
Fourier support in $\bigcup_i\Lex{A_i}$.  If
\[
 \Gamma\setminus\bigcup_i\Lex{A_i}
\]
is finite, then both $\mathcal H^{\exists}_{\mathcal A,\eps}$ and
$\mathcal H^{\forall,\delta}_{\mathcal A,\eps}$ are bounded on
$L^p(\widehat\Gamma)$ for every $1<p<\infty$, uniformly in the coefficients.
In particular, a $4\delta$-separated partition of a word sphere gives a
global $L^p$-unconditional decomposition after adjoining the finite
low-frequency remainder.

The key ingredient of Theorem A is the following modified high-frequency
Cotlar identity, which follows from the hyperbolicity of $X$.  For all
$g,h\in\bigcup_i\Lex{A_i}$ satisfying
$|gh^{-1}|\geq 2m+2\delta$, one has:
\begin{equation}\label{cotlar}
\begin{aligned}
 H^{\exists}(\lambda_g)(H^{\exists})^\circ(\lambda_{h^{-1}})
 ={}&H^{\forall}\big(\lambda_g (H^{\forall})^\circ(\lambda_{h^{-1}})\big)\\
 &+(H^{\forall})^\circ\big(H^{\forall}(\lambda_g)\lambda_{h^{-1}}\big)
 -(H^{\forall})^\circ H^{\forall}(\lambda_g\lambda_{h^{-1}}).
\end{aligned}
\end{equation}
Here,
\[
 H^{\exists}=\mathcal H^{\exists}_{\mathcal A,\eps},
 \qquad
 H^{\forall}=\mathcal H^{\forall,\delta}_{\mathcal A,\eps}.
\]

For $\rho\geq0$ and $m\in\mathbb N$, we say that the Cayley action has
the \emph{$(\rho,m)$-Cotlar property} if, for every finite family
$(A_i)\subset S_m(e)$ with pairwise disjoint branches
$\Lthick{A_i}{\rho}$, every finitely supported coefficient family
$(\eps_i)$ with $|\eps_i|\leq1$, and all
$g,h\in\bigcup_i\Lex{A_i}$ satisfying
$|gh^{-1}|\geq 2m+2\rho$, identity \eqref{cotlar} holds with
\[
 H^{\exists}=\mathcal H^{\exists}_{\mathcal A,\eps},
 \qquad
 H^{\forall}=\mathcal H^{\forall,\rho}_{\mathcal A,\eps}.
\]

\paragraph{Theorem B: Cotlar characterizations of
hyperbolicity (Theorem~\ref{thm:cotlar-characterization}  and Corollary \ref{cor:cayley-cotlar-characterization}).}
 $\Gamma$ is hyperbolic if and only if it has the
$(\rho,m_k)$-Cotlar property   for some $\rho>0$ and an increasing sequence
$m_k\to\infty$ satisfying
\[
 \sup_k\frac{m_{k+1}}{m_k}<\infty.
\]
The same theorem gives equivalent geometric formulations in terms of a
long-product alternative for four chosen geodesics.  Thus the modified
Cotlar identity is not merely a consequence of hyperbolicity: at a
multiplicatively dense collection of scales it detects hyperbolicity.  

In the exact $\rho=0$ setting, if the first-sphere Cotlar identity holds for all distinct relevant pairs, then the Cayley graph is a tree.  Section~5 records a complementary exact, unthickened Cotlar criterion
for general partitions: if
\[
 \Gamma=\bigsqcup_{i\in I_0}\mathcal E_i
\]
has bounded Gromov interface, then, for every  
family $(\eps_i)_{i\in I_0}$ with $|\eps_i|\leq1$, the signed Fourier
multiplier
\[
 \sum_{i\in I_0}\eps_iP_{\mathcal E_i}
\]
is bounded on $L^p(\widehat\Gamma)$ for every $1<p<\infty$
(Theorem~\ref{thm:bounded-interface}).

Exact noncommutative Cotlar identities and geometric sources for them were
studied by Gonz\'alez-P\'erez, Parcet, and Xia \cite{GPX22}; related
identities for Coxeter groups, buildings, and graph products appear in
\cite{LX25,LX26}.  In contrast, the branch identity used here involves a thin/thick pair of
symbols.  The scale-wise converse in Theorem B shows that this modified form is
naturally adapted to general hyperbolic geometry.

The paper is organized as follows.  Section~1 introduces the
branch notation and the relevant Fourier multipliers.  Section~2 proves the
geometric localization and Cotlar identity.  Section~3 derives the local
$L^4$ and global $L^p$ results and records examples.  Section~4 proves the
scale-wise characterization of hyperbolicity.  Section 5 records a complementary bounded-Gromov-interface criterion, its compact-commutator consequence, and several examples.

\section{Notation and branch multipliers}

We will always assume $\Gamma$ is a discrete group  and $\Gamma\curvearrowright X$ is an action by isometries.  We fix a
base point $o\in X$ and write
\[
 |x|=\dist(o,x),\qquad x\in X.
\]
For $r>0$, put
\[
 B_r(a)=\{x\in X:\dist(x,a)<r\},\qquad
 S_r(a)=\{x\in X:\dist(x,a)=r\}.
\]
We adopt the convention
$B_0(a)=S_0(a)=\{a\}$ and $\bar B_r(a)=\{x\in X:\dist(x,a)\leq r\}$.  A geodesic segment joining $x$ and $y$ is denoted by
$[x,y]$.  The Gromov product is
\[
 (y|z)_x=\frac{\dist(x,y)+\dist(x,z)-\dist(y,z)}2.
\]

We impose no global properness assumption on the action.  For $R>0$ put
\begin{equation}
 E_{<R}=\{g\in\Gamma:|g.o|<R\},
 \qquad E_{\geq R}=\Gamma\setminus E_{<R}. \label{ER}
\end{equation}
The low-orbit-frequency hypotheses needed below will be stated explicitly.

For $x\in X$ and $A\subset X$, write
\[
 x\gtrsim A
\]
if at least one geodesic from $o$ to $x$ meets $A$, and write
\[
 x\geq A
\]
if every geodesic from $o$ to $x$ meets $A$.  For a singleton $A=\{a\}$,
we simply write $x\gtrsim a$ and $x\geq a$.  Define
\begin{align*}
 \Lex{A}&=\{g\in\Gamma:g.o\gtrsim A\},\\
 \Lall{A}&=\{g\in\Gamma:g.o\geq A\},\\
 \Lthick{A}{r}&=\{g\in\Gamma:\text{there is }a\in A\text{ such that }
                      g.o\geq B_r(a)\}.
\end{align*}
Also put
\[
 A^r=\{x\in X:\dist(x,A)<r\}=\bigcup_{a\in A}B_r(a).
\]
Note that $\Lthick{A}{r}$ and $\Lall{A^r}$ may not be the same.  In
$\Lthick{A}{r}$, one fixed point $a\in A$ works for all geodesics, whereas
in $\Lall{A^r}$ the intersection point may lie near different points of $A$
on different geodesics.  We always have
\[
 \Lthick{A}{r}\subset \Lall{A^r}.
\]
We say that $X$ is $\delta$-{\it hyperbolic} if every side of every geodesic
triangle lies in the open $\delta$-neighborhood of the union of the other two
sides.
In particular,   any
two geodesics with the same endpoints are $\delta$-close.  So, for $\delta$-hyperbolic $X$, one has 
\[
 \Lall{A}\subset\Lex{A}\subset \Lthick{A}{\delta}\subset \Lall{A^\delta}.
\]

Let $\lambda_g$ denote the left regular representation of $\Gamma$ on
$\ell_2(\Gamma)$.  For a Fourier polynomial
$x=\sum_g c_g\lambda_g$, the canonical trace is given by
\[
 \tau(x)=c_e.
\]
For $1\leq p<\infty$, the associated norm is
\[
 \|x\|_p=\big[\tau(|x|^p)\big]^{1/p},
\]
The resulting Banach space is called the noncommutative $L^p$-space associated with $\Gamma$, which will be denoted by
$L^p(\widehat\Gamma)$.  The theory of $L^p(\hat\Gamma)$ goes back to Segal and Dixmier and has been essential in the recent study of operator approximation properties (\cite{CGPT23},\cite{LS11},\cite{Pa26}); we refer
the reader to \cite{PX03} for background.  When $\Gamma$ is abelian, for
example $\Gamma=\mathbb Z^d$, this space is isometrically isomorphic to the
classical $L^p$-space on the dual group.  The usual duality and complex
interpolation properties remain valid for $L^p(\hat \Gamma)$.  In particular,
$L^p(\widehat\Gamma)^*=L^q(\widehat\Gamma)$ when
$1/p+1/q=1$, and
\[
 [L^p(\widehat\Gamma),L^q(\widehat\Gamma)]_\theta
   =L^r(\widehat\Gamma),
 \qquad
 \frac1r=\frac{1-\theta}{p}+\frac\theta q.
\]
Moreover,
\[
 \|x\|_2=\left(\sum_g|c_g|^2\right)^{1/2}.
\]
Thus, for a function $\varepsilon:\Gamma\to\mathbb C$, the Fourier multiplier
\[
 T_\varepsilon\left(\sum_g c_g\lambda_g\right)
   =\sum_g \varepsilon(g)c_g\lambda_g
\]
is bounded on $L^2(\widehat\Gamma)$ if and only if $\varepsilon$ is bounded.  Its
boundedness on $L^p(\widehat\Gamma)$ for $p\neq2$ is substantially more
subtle and is a central question in harmonic analysis.

  For
$\mathcal E\subset\Gamma$, let $P_{\mathcal E}$ be the orthogonal Fourier
projection onto the closed span of $\{\lambda_g:g\in\mathcal E\}$.

Fix a family $\mathcal A=(A_i)_{i\in I}$.  For $r\geq0$, put
$\mathcal A^r=(A_i^r)_{i\in I}$, with the convention
$\mathcal A^0=\mathcal A$.  For a finitely supported coefficient family
$\eps=(\eps_i)_{i\in I}$ with $|\eps_i|\leq1$, define
\begin{align}
 \mathcal H^{\exists}_{\mathcal A,\eps}
   &=\sum_{i\in I}\eps_iP_{\Lex{A_i}},\label{eq:def-H-exists}\\
 \mathcal H^{\forall,r}_{\mathcal A,\eps}
   &=\sum_{i\in I}\eps_iP_{\Lthick{A_i}{r}}.
   \label{eq:def-H-forall}
\end{align}
Thus geometric thickening is recorded in the subscript,
\[
 \mathcal H^{\exists}_{\mathcal A^r,\eps}
   =\sum_{i\in I}\eps_iP_{\Lex{A_i^r}},
\]
whereas thickening the branch condition is recorded by the superscript
$\forall,r$ in \eqref{eq:def-H-forall}.  Their right-handed companions are
\begin{align*}
 (\mathcal H^{\exists}_{\mathcal A,\eps})^\circ
   &=\sum_{i\in I}\overline{\eps_i}P_{(\Lex{A_i})^{-1}},\\
 (\mathcal H^{\forall,r}_{\mathcal A,\eps})^\circ
   &=\sum_{i\in I}\overline{\eps_i}P_{(\Lthick{A_i}{r})^{-1}}.
\end{align*}
Then
\[
 \big(\mathcal H^{\exists}_{\mathcal A,\eps}(x)\big)^*
   =(\mathcal H^{\exists}_{\mathcal A,\eps})^\circ(x^*),
 \qquad
 \big(\mathcal H^{\forall,r}_{\mathcal A,\eps}(x)\big)^*
   =(\mathcal H^{\forall,r}_{\mathcal A,\eps})^\circ(x^*).
\]

For later use, set
\begin{align}
 \kappa_q(R)
 &=\|P_{E_{<R}}:L^1(\widehat\Gamma)\longrightarrow
                         L^q(\widehat\Gamma)\|,
 \label{eq:kappa}\\
 \eta_q(R)
 &=\|P_{E_{\geq R}}:L^q(\widehat\Gamma)\longrightarrow
                         L^q(\widehat\Gamma)\|.
 \label{eq:eta}
\end{align}
If $E_{<R}$ is finite, then both quantities are finite.  

\section{Geometry and the Cotlar property}

For   Sections~2 and 3, we assume that $X$ is
$\delta$-hyperbolic. 
We begin with an elementary form of thin-triangle geometry.

\begin{lemma}[Fellow traveling]\label{lem:fellow-travel}
 
 For every $x,y,z\in X$, and
      geodesic $\alpha=[x,y]$, $\beta=[x,z]$, suppose   $a\in\alpha$ satisfies
      \begin{equation}\label{eq:gromov-threshold}
       \dist(x,a)+\delta\leq (y|z)_x,
      \end{equation}
      then
      $$\beta\cap B_\delta(a)\neq\varnothing.$$ 
\end{lemma}

\begin{proof}
   Complete $\alpha$ and
$\beta$ to a geodesic triangle by choosing a geodesic $\gamma=[y,z]$.
The hyperbolicity of $X$  places $a$ within distance $\delta$ of $\beta\cup\gamma$.  Suppose $\dist(\beta,a)\geq \delta$, then there exists   $b\in\gamma$ and $\dist(a,b)<\delta$.  We obtain
\begin{align*}
 \dist(x,z)
 &\leq \dist(x,a)+\dist(a,b)+\dist(b,z)\\
 &= \dist(x,a)+\dist(a,b)+\dist(y,z)-\dist(b,y)\\
 &\leq \dist(x,a)+\dist(a,b)+\dist(y,z)-\dist(y,a)+\dist(a,b)\\
 &= \dist(x,a)+\dist(a,b)+\dist(y,z)-(\dist(x,y)-\dist(x,a))+\dist(a,b)\\ 
 &=2\dist(x,a)+\dist(y,z)-\dist(x,y)+2\dist(a,b).
\end{align*}
Thus $(y|z)_x\leq \dist(x,a)+\dist(a,b)<\dist(x,a)+\delta$, contradicting
\eqref{eq:gromov-threshold}.  Hence $\dist(\beta,a)<\delta$.

\end{proof}

\medskip
 
\begin{lemma}[Branch-propagation]\label{lem:branch-propagation}
For every $g,h\in\Gamma$ and
      $a\in X$, the conditions
      \begin{equation}\label{eq:branch-propagation-hyp}
       g.o\gtrsim a,
       \qquad
       \frac{|gh.o|+|g.o|-|h.o|}{2}\geq |a|+\delta
      \end{equation}
      imply
      \begin{equation}\label{eq:branch-propagation-conclusion}
       gh.o\geq B_\delta(a).
      \end{equation} 
 
\end{lemma}

\begin{proof}
 Choose a geodesic
$\alpha=[o,g.o]$ passing through $a$ and an arbitrary geodesic
$\beta=[o,gh.o]$.  Since the action is isometric,
\[
 \dist(g.o,gh.o)=\dist(o,h.o)=|h.o|.
\]
Thus \eqref{eq:branch-propagation-hyp} says precisely that
\[
 \dist(o,a)+\delta\leq (g.o\,|\,gh.o)_o.
\]
Lemma~\ref{lem:fellow-travel} therefore shows that $\beta$ meets
$B_\delta(a)$.  Since $\beta$ was arbitrary,
\eqref{eq:branch-propagation-conclusion} follows.

\end{proof}

\begin{lemma}[Localization of long products]\label{lem:localization}
Let $m>0$ and let $(A_i)_{i\in I}\subset\overline B_m(o)$.  If
$g\in\Lex{A_i}$, $h\in\Lex{A_j}$, and $|gh^{-1}.o|\geq 2m+2\delta $, 
then
\[
 gh^{-1}\in\Lthick{A_i}{\delta}
 \qquad\text{or}\qquad
 hg^{-1}\in\Lthick{A_j}{\delta}.
\]
Consequently, for any $i,j\in I$,
\begin{equation}\label{eq:localization}
 E_{\geq 2m+2\delta}
 \cap\big(\Lex{A_i}(\Lex{A_j})^{-1}\big)
 \subset \Lthick{A_i}{\delta}\cup (\Lthick{A_j}{\delta})^{-1}.
\end{equation}
\end{lemma}

\begin{proof}
Let $g\in\Lex{A_i}$ and $h\in\Lex{A_j}$, and put $q=gh^{-1}$.  Assume
$|q.o|\geq2m+2\delta$.  If $|g.o|\geq|h.o|$, choose $a\in A_i$ on a
geodesic from $o$ to $g.o$.  Then
\[
 \frac{|q.o|+|g.o|-|h.o|}{2}\geq m+\delta\geq |a|+\delta.
\]
Lemma~\ref{lem:branch-propagation}, applied to $(g,h^{-1})$, gives
$q\in \Lthick{A_i}{\delta}$.  If $|h.o|>|g.o|$, apply the same argument to
$(h,g^{-1})$ to obtain $q^{-1}\in \Lthick{A_j}{\delta}$.
\end{proof}

\begin{proposition}[Cotlar identity for separated branches]\label{lem:cotlar}
Let $m>0$, and let $(A_i)_{i\in I}\subset S_m(o)$ satisfy
\begin{equation}\label{eq:thick-disjoint}
 \Lthick{A_i}{\delta}\cap \Lthick{A_j}{\delta}=\varnothing,
 \qquad i\neq j.
\end{equation}
For all $g,h\in\bigcup_i\Lex{A_i}$ satisfying
$gh^{-1}\in E_{\geq2m+2\delta}$, the following branch Cotlar identity
holds:
\begin{equation}\label{eq:cotlar}
\begin{aligned}
 H^{\exists}(\lambda_g)(H^{\exists})^\circ(\lambda_{h^{-1}})
 ={}&H^{\forall}\big(\lambda_g (H^{\forall})^\circ(\lambda_{h^{-1}})\big)
 +(H^{\forall})^\circ\big(H^{\forall}(\lambda_g)\lambda_{h^{-1}}\big)
 -(H^{\forall})^\circ H^{\forall}(\lambda_g\lambda_{h^{-1}}),\\
 ={}&H^{\forall}\big(\lambda_g (H^{\exists})^\circ(\lambda_{h^{-1}})\big)
 +(H^{\forall})^\circ\big(H^{\exists}(\lambda_g)\lambda_{h^{-1}}\big)
 -(H^{\forall})^\circ H^{\forall}(\lambda_g\lambda_{h^{-1}}).
\end{aligned}
\end{equation}
Here,
\[
 H^{\exists}=\mathcal H^{\exists}_{\mathcal A,\eps},
 \qquad
 H^{\forall}=\mathcal H^{\forall,\delta}_{\mathcal A,\eps}.
\]
\end{proposition}

\begin{proof} Note that $\Lex{A_j}\subset \Lthick{A_j}{\delta}$ and $H^{\exists}$ and $H^{\forall}$ coincide on $\bigcup_i\Lex{A_i}$. So the two lines on the right hand side of \eqref{eq:cotlar} are equal.
  Put $q=gh^{-1}$ and choose
$i,j$ with $g\in\Lex{A_i}$ and $h\in\Lex{A_j}$.  By
Lemma~\ref{lem:localization}, either
$q\in \Lthick{A_i}{\delta}$ or $q^{-1}\in \Lthick{A_j}{\delta}$.
 
If $q\in \Lthick{A_i}{\delta}$, then $  H^{\exists}(\lambda_g)\lambda_{h^{-1}}
 = H^{\forall}(\lambda_g\lambda_{h^{-1}})=\varepsilon_i\lambda_q$ and \begin{eqnarray*}   H^{\forall}\big(\lambda_g (H^{\exists})^\circ(\lambda_{h^{-1}})\big)&= &H^{\exists}(\lambda_g)(H^{\exists})^\circ(\lambda_{h^{-1}})=\varepsilon_i\bar\varepsilon_j\lambda_q\\ 
(H^{\forall})^\circ\big(H^{\exists}(\lambda_g)\lambda_{h^{-1}}\big)
 &=&(H^{\forall})^\circ H^{\forall}(\lambda_g\lambda_{h^{-1}})=(H^{\forall})^\circ(\varepsilon_i \lambda_q).
 \end{eqnarray*}
Hence \eqref{eq:cotlar}.

  If
$q^{-1}\in \Lthick{A_j}{\delta}$, then $   \lambda_g (H^{\exists})^\circ\lambda_{h^{-1}}
 =( H^{\forall})^\circ(\lambda_g\lambda_{h^{-1}})=\bar\varepsilon_j\lambda_q$ and 
 \begin{eqnarray*}   H^{\forall}\big(\lambda_g (H^{\exists})^\circ(\lambda_{h^{-1}})\big)&= &(H^{\forall})^\circ H^{\forall}(\lambda_g\lambda_{h^{-1}})\\ 
(H^{\forall})^\circ\big(H^{\exists}(\lambda_g)\lambda_{h^{-1}}\big)
 &=&H^{\exists}(\lambda_g)(H^{\exists})^\circ(\lambda_{h^{-1}}).
 \end{eqnarray*} Hence \eqref{eq:cotlar}.
\end{proof}

\begin{remark}
We will show in Section~4 that, for connected vertex-transitive graphs,
the scale-wise branch Cotlar property associated with
Proposition~\ref{lem:cotlar} characterizes hyperbolicity.
\end{remark} 

\section{\texorpdfstring{$L^p$}{Lp}-unconditionality}

\begin{theorem}[Main theorem]\label{thm:main}
Let $m>0$, let $(A_i)_{i\in I}\subset S_m(o)$ satisfy
\eqref{eq:thick-disjoint}, and put $R=2m+2\delta$.  Assume that the set $E_{<R}$ defined in (\ref{ER}) is
finite.  Then there is a constant $C_4$, depending only on
$\kappa_2(R)$, such that
\begin{equation}\label{eq:L4-main}
 \|\mathcal H^{\exists}_{\mathcal A,\eps}x\|_4
 =\|\mathcal H^{\forall,\delta}_{\mathcal A,\eps}x\|_4
 \leq C_4\|x\|_4
\end{equation}
for every Fourier polynomial $x$ supported in
$\bigcup_i\Lex{A_i}$ and every $|\eps_i|\leq1$.  Consequently, the
subspaces
\[
 X_i^4=\overline{\operatorname{span}}^{\,L^4}
       \{\lambda_g:g\in\Lex{A_i}\}
\]
form an $L^4$-unconditional family on their closed linear span.  All
estimates are uniform over finite subfamilies and finitely supported
coefficient choices.

Assume additionally that
\begin{equation}\label{eq:finite-complement}
 F=\Gamma\setminus\bigcup_i\Lex{A_i}
\end{equation}
is finite.  Then, for every $1<p<\infty$, there is a constant $C_p$ such
that
\begin{equation}\label{eq:Lp-main}
 \|\mathcal H^{\exists}_{\mathcal A,\eps}x\|_p
 +\|\mathcal H^{\forall,\delta}_{\mathcal A,\eps}x\|_p
 \leq C_p\|x\|_p,
 \qquad x\in L^p(\widehat\Gamma).
\end{equation}
Thus, after adjoining the finite-dimensional subspace generated by
$\{\lambda_g:g\in F\}$, the branch subspaces form an unconditional
decomposition of $L^p(\widehat\Gamma)$.
\end{theorem}

\begin{proof}
Write
\[
 H^{\exists}=\mathcal H^{\exists}_{\mathcal A,\eps},
 \qquad
 H^{\forall}=\mathcal H^{\forall,\delta}_{\mathcal A,\eps}.
\]
On the Fourier support of $x$ we have
$H^{\forall}x=H^{\exists}x$.  Since the sets
$\Lthick{A_i}{\delta}$ are pairwise disjoint and
$\Lex{A_i}\subset \Lthick{A_i}{\delta}$, the four branch multipliers
$H^{\exists},(H^{\exists})^\circ,H^{\forall},(H^{\forall})^\circ$ and the
projection $P_{E_{\geq R}}$ are contractions on $L^2$.

We decompose
\[
 |(H^{\exists}x)^*|^2
 =P_{E_{<R}}|(H^{\exists}x)^*|^2
  +P_{E_{\geq R}}\big(H^{\exists}x\,(H^{\exists})^\circ(x^*)\big).
\]
For the low-frequency term,
\begin{align*}
 \|P_{E_{<R}}|(H^{\exists}x)^*|^2\|_2
 &\leq \kappa_2(R)\,\|| (H^{\exists}x)^*|^2\|_1\\
 &=\kappa_2(R)\|H^{\exists}x\|_2^2
 \leq\kappa_2(R)\|x\|_4^2.
\end{align*}
After summing Proposition~\ref{lem:cotlar} over the Fourier coefficients of
$x$, the high-frequency term equals
\[
 P_{E_{\geq R}}\Big(
 H^{\forall}\big(x(H^{\forall})^\circ(x^*)\big)
 +(H^{\forall})^\circ\big(H^{\forall}(x)x^*\big)
 -(H^{\forall})^\circ H^{\forall}(xx^*)
 \Big).
\]
Hence, by the $L^2$-contractivity of the three multipliers and
H\"older's inequality,
\begin{align*}
 &\left\|P_{E_{\geq R}}
    \big(H^{\exists}x\,(H^{\exists})^\circ(x^*)\big)\right\|_2\\
 &\qquad\leq
 \|x(H^{\forall})^\circ(x^*)\|_2
 +\|H^{\forall}(x)x^*\|_2+\|xx^*\|_2\\
 &\qquad\leq 2\|x\|_4\|H^{\forall}x\|_4+\|x\|_4^2\\
 &\qquad=2\|x\|_4\|H^{\exists}x\|_4+\|x\|_4^2.
\end{align*}
Therefore
\[
 \|H^{\exists}x\|_4^2
 \leq (\kappa_2(R)+1)\|x\|_4^2
      +2\|x\|_4\|H^{\exists}x\|_4.
\]
Solving this quadratic inequality gives, for example,
\[
 \|H^{\exists}x\|_4
 \leq \big(1+\sqrt{\kappa_2(R)+2}\big)\|x\|_4.
\]
This proves \eqref{eq:L4-main}.  The two-sided unconditional estimate for
unimodular coefficients follows by applying the same estimate to
$\sum_i\eps_i x_i$ with coefficients $\overline{\eps_i}$.

Now assume \eqref{eq:finite-complement}.  Let
$P_0=P_{\cup_i\Lex{A_i}}=I-P_F$.  Since $F$ is finite, $P_0$ is
bounded on every $L^p$.  Applying the local $L^4$ estimate to $P_0x$
shows that $H^{\exists}$ is bounded on all of $L^4$.  Moreover,
$D=H^{\forall}-H^{\exists}$ has Fourier support contained in $F$, so $D$
is finite rank and $H^{\forall}$ is also bounded on $L^4$.

We spell out the dyadic induction.  Suppose that $H^{\exists}$ and
$H^{\forall}$ are bounded on $L^q$, where $q=2^k\geq4$, uniformly over
the coefficient choices, and put
\[
 M_q=\max\{\|H^{\exists}\|_{L^q\to L^q},
             \|H^{\forall}\|_{L^q\to L^q}\}.
\]
The right-handed multipliers have the same $L^q$-norms.  For a Fourier
polynomial $x$, put $x_0=P_0x$.  Then
$H^{\exists}x=H^{\exists}x_0=H^{\forall}x_0$.  Repeating the preceding
decomposition in $L^q$, and retaining the norm of the high-frequency
projection, gives
\begin{align}
 \|H^{\exists}x\|_{2q}^2
 &\leq \kappa_q(R)\|x_0\|_{2q}^2
 \notag\\
 &\quad+\eta_q(R)\Big(
       2M_q\|x_0\|_{2q}\|H^{\exists}x\|_{2q}
       +M_q^2\|x_0\|_{2q}^2\Big).
 \label{eq:dyadic-induction}
\end{align}
For example,
\begin{align*}
 \|H^{\forall}\big(x_0(H^{\forall})^\circ(x_0^*)\big)\|_q
 &\leq M_q\|x_0\|_{2q}
       \|(H^{\forall})^\circ(x_0^*)\|_{2q}\\
 &=M_q\|x_0\|_{2q}\|H^{\exists}x\|_{2q},
\end{align*}
and the second term is estimated in the same way, while the third is
bounded by $M_q^2\|x_0\|_{2q}^2$.  Since $P_0$ is bounded on
$L^{2q}$, the quadratic inequality \eqref{eq:dyadic-induction} proves
the $L^{2q}$-boundedness of $H^{\exists}$.  The identity
$H^{\forall}=H^{\exists}+D$ then gives the $L^{2q}$-boundedness of
$H^{\forall}$.  Starting with $q=4$, this proves boundedness at every
dyadic exponent $2^k$, $k\geq2$.  Interpolation between $L^2$ and the
dyadic spaces gives \eqref{eq:Lp-main} for $2\leq p<\infty$.  Finally,
duality gives the range $1<p\leq2$.
\end{proof}

\begin{remark}\label{rem:main-low-frequency}
The finiteness of $E_{<R}$ in Theorem~\ref{thm:main} is a convenient
sufficient hypothesis and is automatic when the orbit map is metrically proper, that is, when every set $E_{<R}$ is finite.  We avoid
imposing properness in order to retain the generality of the metric-action
framework.  For the $L^4$ conclusion, the finiteness of $E_{<R}$ may be
replaced by $\kappa_2(R)<\infty$.  For the dyadic induction leading to the
global $L^p$ conclusion, it is enough that $\kappa_q(R)$ and $\eta_q(R)$ be
finite at the dyadic exponents that occur.
\end{remark}

\begin{corollary}[Separated subsets of a sphere]\label{cor:sphere}
Let $m>0$ and let $(A_i)_{i\in I}\subset S_m(o)$ satisfy
\begin{equation}\label{eq:sphere-separation}
 \dist(A_i,A_j)\geq4\delta,
 \qquad i\neq j.
\end{equation}
Assume that $E_{<2m+2\delta}$ is finite.  Then \eqref{eq:L4-main}
 holds.
\end{corollary}

\begin{proof}
Suppose that $s\in \Lthick{A_i}{\delta}\cap \Lthick{A_j}{\delta}$.  Choose one
geodesic $\sigma$ from $o$ to $s.o$.  There are $a\in A_i$, $b\in A_j$,
and $u,v\in\sigma$ such that
\[
 \dist(a,u)<\delta,
 \qquad
 \dist(b,v)<\delta.
\]
Because $u$ and $v$ lie on the same geodesic issuing from $o$,
\[
 \dist(u,v)=\big||u|-|v|\big|<2\delta,
\]
where we used $|a|=|b|=m$.  Thus $\dist(a,b)<4\delta$, contradicting
\eqref{eq:sphere-separation}.  Hence the thick branches are pairwise
disjoint, and Theorem~\ref{thm:main} applies.
\end{proof}

\begin{corollary}[A separated partition of a sphere]\label{cor:sphere-partition}
Let $m>0$.  Suppose that
\[
 S_m(o)=\bigsqcup_{i\in I}A_i,
 \qquad
 \dist(A_i,A_j)\geq4\delta\quad(i\neq j),
\]
and that $E_{<2m+2\delta}$ is finite.
Then the corresponding branch subspaces are $L^p$-unconditional for every
$1<p<\infty$, after adjoining the finite-dimensional space generated by
$E_{<m}$.  Equivalently, the signed branch multiplier is bounded on
$L^p(\widehat\Gamma)$, uniformly over all finitely supported coefficient choices
$|\eps_i|\leq1$.
\end{corollary}

\begin{proof}
Corollary~\ref{cor:sphere} gives the separation hypothesis of
Theorem~\ref{thm:main}.  Moreover, every geodesic from $o$ to $g.o$ with
$|g.o|\geq m$ meets $S_m(o)$, so
\[
 E_{\geq m}\subset\bigcup_i\Lex{A_i}.
\]
Thus the complement of the union of the branches is contained in the finite
set $E_{<m}$, and the global part of Theorem~\ref{thm:main} applies.
\end{proof}

\subsection{Examples}\label{subsec:examples}

We record several concrete instances of the separated-branch construction.
The first two give global $L^p$-unconditional decompositions, while the last
two illustrate the local $L^4$ conclusion of Theorem~\ref{thm:main}.

\paragraph{Tree actions and prefix branches.}
Let $T$ be a locally finite tree, and let $\Gamma\curvearrowright T$ act by
graph automorphisms.  Assume that the stabilizer $\Gamma_o$ of a vertex
$o\in T$ is finite.  Then $E_{<R}$ is finite for every $R>0$: the vertex
ball $B_R(o)$ is finite, and every nonempty fiber of the orbit map
$g\mapsto g.o$ is a coset of $\Gamma_o$.  Fix $m\in\N$ and any partition
\[
 S_m(o)=\bigsqcup_{i\in I}A_i.
\]
By the radius-zero specialization stated after
Lemma~\ref{lem:fellow-travel}, Corollary~\ref{cor:sphere-partition} gives an $L^p$-unconditional
decomposition for every $1<p<\infty$, after adjoining the finite-dimensional
space generated by $E_{<m}$.

For the standard Cayley tree of $\mathbb F_n$, taking the $A_i$ to be
singletons gives the decomposition according to prefixes of length $m$:
$\Lex{\{a\}}$ consists of the reduced words whose first $m$ letters form
$a\in S_m(e)$.  The case $m=1$ is the usual first-letter decomposition.
The same argument applies, for example, to
$\mathbb Z_2*\mathbb Z_3$ acting on its Bass--Serre tree.

\begin{corollary}[End sectors]\label{cor:end-sectors}
Let $\Gamma$ be a finitely generated group whose Cayley graph
$X=\operatorname{Cay}(\Gamma,S)$ is $\delta$-hyperbolic.  Choose integers
$0\leq r<m$ satisfying
\[
 m-r\geq2\delta.
\]
Let $(C_i)_{i\in I}$ be the connected components of
$X\setminus\overline B_r(e)$ that meet $S_m(e)$, and put
\[
 A_i=C_i\cap S_m(e).
\]
Then the branch subspaces associated with $(A_i)_{i\in I}$ form an
$L^p$-unconditional decomposition of $L^p(\widehat\Gamma)$, after adjoining
the finite-dimensional space generated by $E_{<m}$, for every
$1<p<\infty$.
\end{corollary}

\begin{proof}
The sets $(A_i)_{i\in I}$ form a partition of $S_m(e)$.  If
$a\in A_i$ and $b\in A_j$ with $i\neq j$, every path from $a$ to $b$ meets
$\overline B_r(e)$.  Consequently,
\[
 \dist(a,b)\geq \dist(a,\overline B_r(e))
                  +\dist(b,\overline B_r(e))
              =2(m-r)\geq4\delta.
\]
Since word balls in a finitely generated group are finite,
Corollary~\ref{cor:sphere-partition} applies.
\end{proof}

If $\Gamma$ has more than one end, one may choose $r$ so that
$X\setminus\overline B_r(e)$ has at least two unbounded components.  For all
sufficiently large $m$, Corollary~\ref{cor:end-sectors} then gives a
nontrivial unconditional decomposition indexed by the corresponding end
sectors.

\paragraph{Separated shadows in a hyperbolic Cayley graph.}
Let $\Gamma$ be any finitely generated $\delta$-hyperbolic group and choose
a finite set
\[
 \{a_1,\ldots,a_N\}\subset S_m(e),
 \qquad
 \dist(a_i,a_j)\geq4\delta\quad(i\neq j).
\]
Taking $A_i=\{a_i\}$ in Corollary~\ref{cor:sphere} gives
\[
 \left\|\sum_{i=1}^N\eps_iP_{\Lex{\{a_i\}}}x\right\|_4
 \leq C_4\|x\|_4
\]
for every Fourier polynomial $x$ supported in
$\bigcup_i\Lex{\{a_i\}}$, uniformly over $|\eps_i|\leq1$.  Here
$\Lex{\{a_i\}}$ is the geodesic shadow of $a_i$.  This construction does
not require a decomposition into ends and therefore also gives nontrivial
local examples for one-ended hyperbolic groups.

\paragraph{Angular sectors in real hyperbolic space.} Let \(\mathbb H^d\) denote the \(d\)-dimensional real hyperbolic space of constant sectional curvature $-1$.
Let $\Gamma<\operatorname{Isom}(\mathbb H^d)$ be a discrete group acting
properly, fix $o\in\mathbb H^d$, and choose a hyperbolicity constant
$\delta$ for $\mathbb H^d$.  Let
$A_1,\ldots,A_N\subset S_m(o)$ be spherical caps, or more general subsets,
satisfying
\[
 \dist(A_i,A_j)\geq4\delta\quad(i\neq j).
\]
The orbit ball $E_{<2m+2\delta}$ is finite, so
Corollary~\ref{cor:sphere} gives an $L^4$-unconditional family on the closed
span of the corresponding branches.  Since $\mathbb H^d$ is uniquely
geodesic, $\Lex{A_i}$ consists exactly of those $g\in\Gamma$ for which the
radial geodesic from $o$ to $g.o$ crosses $A_i$.  A concrete example is the action on $\mathbb H^d$ of the fundamental
group of a closed hyperbolic manifold.
This is typically a local $L^4$ example rather than a global decomposition,
since separated caps need not cover the whole metric sphere.

\section{Cotlar characterizations of hyperbolicity}
\label{sec:cotlar-characterization}

Throughout this section, $X$ is a connected graph equipped with its path
metric and viewed as a metric realization with every edge of length one.
The group $\Gamma$ acts transitively on the vertex set $V(X)$ by graph
automorphisms, and $o\in V(X)$ is fixed.  Here and below, transitivity means
vertex-transitivity; no properness assumption is imposed.  Since every point
of the metric realization is uniformly close to a vertex, replacing points
by vertices does not change any asymptotic cone.

We show that requiring the branch Cotlar identity of Proposition \ref{lem:cotlar} at every scale also characterizes hyperbolicity.

\begin{definition}[Cotlar property for actions]\label{def:cotlar-property}
Let $\rho\geq0$ and $m\in\N$.  We say that the action of $\Gamma$ on
$X$ has the \emph{$(\rho,m)$-Cotlar property} if, for every finite family
$(A_i)_{i\in I}$ of subsets of
$S_m(o)$ satisfying
\[
 \Lthick{A_i}{\rho}\cap \Lthick{A_j}{\rho}=\varnothing,
 \qquad i\neq j,
\]
every choice $|\eps_i|\leq1$, and every
$g,h\in\bigcup_i\Lex{A_i}$ with
\[
 |gh^{-1}.o|\geq2m+2\rho,
\]
the identity \eqref{eq:cotlar} holds with
$H^{\exists}=\mathcal H^{\exists}_{\mathcal A,\eps}$ and
$H^{\forall}=\mathcal H^{\forall,\rho}_{\mathcal A,\eps}$, where the latter
multiplier is defined by \eqref{eq:def-H-forall} with $r=\rho$.  For
$\rho=0$, the convention
$B_0(a)=\{a\}$ is in force.  We say that the action has the
\emph{$\rho$-Cotlar property} if it has the $(\rho,m)$-Cotlar property for
every $m\in\N$.
\end{definition}

Before proving the general characterization, we record the exact tree
argument in the Cayley-graph case.  It explains the scalar factorization that
occurs for free groups, but it is not needed in the positive-radius converse.

\begin{proposition}[The exact Cayley-tree case]
\label{prop:exact-tree}
Let $X=\operatorname{Cay}(\Gamma,S)$ with base point $e$.  Assume that,
for every finite set $F\subset S_1(e)$, the exact $\rho=0$ version of
the Cotlar identity \eqref{eq:cotlar} holds for the singleton family
\[
 (A_a)_{a\in F},\qquad A_a=\{a\},
\]
for every coefficient choice $(\eps_a)_{a\in F}$ with $|\eps_a|\leq1$
and all distinct
\[
 g,h\in\bigcup_{a\in F}\mathcal L^{\exists}_{\{a\}}.
\]
Then $X$ is a tree.
\end{proposition}

\begin{proof}

For $c\in S_1(e)$ and $s\in\Gamma$, put
\[
 u_c(s)=\mathbf 1_{\mathcal L^{\exists}_{\{c\}}}(s),
 \qquad
 v_c(s)=\mathbf 1_{\mathcal L^{\forall}_{\{c\}}}(s).
\]Suppose, toward a contradiction, that $X$ contains a simple cycle.
Choose a simple cycle $C$ of minimal length.  Choose vertices $x,y,z$
in cyclic order on $C$ so that each of the three arcs of $C$ between
consecutive vertices has length at most $|C|/2$.  Each of these arcs is
geodesic: otherwise, a shorter path between the same endpoints,
together with the corresponding arc, would contain a simple cycle
shorter than $C$.

Put
\[
 h=x^{-1}y,
 \qquad
 g=z^{-1}y,
 \qquad
 q=gh^{-1}=z^{-1}x.
\]
Let $a\in S_1(e)$ be the first vertex of the geodesic obtained by
translating the arc from $z$ to $y$ by $z^{-1}$, and let
$b\in S_1(e)$ be the first vertex of the geodesic obtained by
translating the arc from $x$ to $y$ by $x^{-1}$.  Then
\[
 g\in\Lex{\{a\}},
 \qquad
 h\in\Lex{\{b\}},
\]
so that $u_a(g)=u_b(h)=1$.

The translated arc from $z$ to $x$ is a geodesic from $e$ to
$q=z^{-1}x$.  At $z$, the arcs from $z$ to $y$ and from $z$ to $x$
leave along distinct edges of the simple cycle.  Hence this geodesic
has first vertex different from $a$, and therefore it does not meet
$a$ at all.  Consequently,
\[
 q\notin\Lall{\{a\}},
 \qquad
 v_a(q)=0.
\]
Similarly, the translated arc from $x$ to $z$ is a geodesic from $e$
to $q^{-1}=x^{-1}z$ which avoids $b$.  Thus
\[
 q^{-1}\notin\Lall{\{b\}},
 \qquad
 v_b(q^{-1})=0.
\]

Apply the assumed identity with $F=\{a,b\}$ if $a\neq b$, and with
$F=\{a\}$ if $a=b$. Evaluate the exact Cotlar identity on $\lambda_g$ and
$\lambda_{h^{-1}}$ and extract the coefficient of
$\eps_a\overline{\eps_b}$.  If $a\neq b$, this coefficient may be
extracted by integrating over the coefficient torus; if $a=b$, set all
other coefficients equal to zero.  We obtain
\begin{equation}\label{eq:exact-tree-coefficient}
 u_a(g)u_b(h)
 =v_a(q)v_b(h)+v_a(g)v_b(q^{-1})-v_a(q)v_b(q^{-1}).
\end{equation}
The left-hand side of \eqref{eq:exact-tree-coefficient} equals $1$,
whereas the right-hand side equals $0$ because
$v_a(q)=v_b(q^{-1})=0$.  This contradiction shows that $X$ contains no
simple cycle.  Since a Cayley graph is connected, $X$ is a tree.
\end{proof}

We recall the asymptotic-cone terminology which will be used in this section.  Fix a
nonprincipal ultrafilter $\omega$, basepoints $(x_n)$, and scaling constants
$r_n\to\infty$.  The asymptotic cone
\[
 \Cone_\omega(X,(x_n),(r_n))
\]
is the metric ultralimit of the pointed spaces $(X,r_n^{-1}d,x_n)$.  A
geodesic, or a geodesic triangle, obtained as the ultralimit of
coordinatewise geodesics is called \emph{internal}.  For an arbitrary
geodesic metric space, hyperbolicity is equivalent to every asymptotic cone
being an $\mathbb R$-tree; see Frigerio--Sisto \cite{FS09} for this precise
generality and a proof.

\begin{lemma}[Fixed-scaling asymptotic-cone criterion]
\label{lem:fixed-scaling}
A geodesic metric space $X$ is hyperbolic if and only if, for every
nonprincipal ultrafilter $\omega$ and every basepoint sequence $(x_n)$, the
asymptotic cone
\[
 \Cone_\omega(X,(x_n),(n))
\]
is an $\mathbb R$-tree.
\end{lemma}

\begin{proof}
The forward implication follows from the standard asymptotic-cone
characterization recalled above.  The converse follows from \cite{DK18} Proposition 11.167.
\end{proof}

\begin{theorem}[Cotlar and long-product characterizations of hyperbolicity]
\label{thm:cotlar-characterization}
Suppose that the group $\Gamma$ acts transitively on the vertex set $V(X)$
by graph automorphisms.  Then the following are equivalent:
\begin{enumerate}[label=\textup{(\roman*)}]
\item $X$ is hyperbolic.

\item The action has the $\rho$-Cotlar property for some $\rho>0$.

\item There exist $\rho>0$ and a strictly increasing sequence of positive
      integers $(m_k)_{k\geq1}$ such that
      \begin{equation}\label{eq:bounded-ratio-scales}
       m_k\longrightarrow\infty,
       \qquad
       C_0:=\sup_{k\geq1}\frac{m_{k+1}}{m_k}<\infty,
      \end{equation}
      and the action has the $(\rho,m_k)$-Cotlar property for every $k$.

\item There exist $\rho>0$ and a sequence $(m_k)_{k\geq1}$ satisfying
      \eqref{eq:bounded-ratio-scales} with the following property.  For every
      $k\geq1$, every $h_1,h_2\in\Gamma$, and every
      \[
       A_1,A_2\subset   S_{m_k}(o)
      \]
      satisfying
      \[
       A_1=A_2
       \qquad\text{or}\qquad
       \dist(A_1,A_2)\geq4\rho,
      \]
      put $q=h_1h_2^{-1}$.  Given geodesics
      \[
      \begin{aligned}
       \gamma_1&=[o,h_1.o], & \gamma_2&=[o,h_2.o],\\
       \gamma_3&=[o,q.o],   & \gamma_4&=[o,q^{-1}.o],
      \end{aligned}
      \]
      if
      \[
       A_1\cap\gamma_1\neq\varnothing,
       \qquad
       A_2\cap\gamma_2\neq\varnothing,
       \qquad
       |q.o|\geq2m_k+2\rho,
      \]
      then
      \begin{equation}\label{eq:set-long-product-alternative}
       \gamma_3\cap A_1^\rho\neq\varnothing
       \qquad\text{or}\qquad
       \gamma_4\cap A_2^\rho\neq\varnothing.
      \end{equation}

\item There exist $\rho>0$ and a sequence $(m_k)_{k\geq1}$ satisfying
      \eqref{eq:bounded-ratio-scales} with the following property.  For every
      $k\geq1$, every $h_1,h_2\in\Gamma$, every
      \[
       a_1,a_2\in   S_{m_k}(o),
      \]
      and geodesics $\gamma_1,\ldots,\gamma_4$ as in \textup{(iv)}, if
      \[
       a_1\in\gamma_1,
       \qquad
       a_2\in\gamma_2,
       \qquad
       |h_1h_2^{-1}.o|\geq2m_k+2\rho,
      \]
      then
      \begin{equation}\label{eq:point-long-product-alternative}
       \gamma_3\cap B_\rho(a_1)\neq\varnothing
       \qquad\text{or}\qquad
       \gamma_4\cap B_\rho(a_2)\neq\varnothing.
      \end{equation}
\end{enumerate}
\end{theorem}

\begin{proof}
\textup{(i)}$\Rightarrow$\textup{(ii)}.
Choose a hyperbolicity constant $\rho>0$ for $X$.  Lemmas
\ref{lem:branch-propagation} and \ref{lem:localization}, together with
Proposition~\ref{lem:cotlar}, show that the action has the $\rho$-Cotlar
property.

\smallskip
\noindent
\textup{(ii)}$\Rightarrow$\textup{(iii)}.
Take $m_k=k$.  Then \eqref{eq:bounded-ratio-scales} holds, and the conclusion
is immediate from the definition of the $\rho$-Cotlar property.

\smallskip
\noindent
\textup{(iii)}$\Rightarrow$\textup{(iv)}.
We first record the geometric consequence of the Cotlar identity at a
permitted scale.  Let $(A_i)$ be an admissible family in
$S_{m_k}(o)$ and put
\[
 u_i(s)=\mathbf1_{\Lex{A_i}}(s),
 \qquad
 v_i(s)=\mathbf1_{\Lthick{A_i}{\rho}}(s).
\]
If $g\in\Lex{A_i}$, $h\in\Lex{A_j}$, and $q=gh^{-1}$ satisfies
$|q.o|\geq2m_k+2\rho$, evaluating \eqref{eq:cotlar} and extracting the
$\eps_i\overline{\eps_j}$ coefficient gives
\begin{equation}\label{eq:coefficient-cotlar}
 u_i(g)u_j(h)
 =v_i(q)v_j(h)+v_i(g)v_j(q^{-1})-v_i(q)v_j(q^{-1}).
\end{equation}
For $i\neq j$, the coefficient is extracted by integrating over the
coefficient torus; for $i=j$, one first sets all other coefficients equal to
zero.  Since the left-hand side equals $1$, it is impossible that both
$v_i(q)=0$ and $v_j(q^{-1})=0$.  Hence
\begin{equation}\label{eq:geometric-cotlar}
 q\in\Lthick{A_i}{\rho}
 \qquad\text{or}\qquad
 q^{-1}\in\Lthick{A_j}{\rho}.
\end{equation}

Now fix the data in \textup{(iv)}.  If $A_1=A_2$, use the one-member family
$(A_1)$.  If $\dist(A_1,A_2)\geq4\rho$, then the two-member family
$(A_1,A_2)$ is admissible.  Indeed, if an element belonged to both thick
branches, there would be $a_i\in A_i$ and points $x_i$ on one geodesic from
$o$ such that $\dist(a_i,x_i)<\rho$.  Since $|a_1|=|a_2|=m_k$,
\[
 \dist(x_1,x_2)=\big||x_1|-|x_2|\big|<2\rho,
\]
and therefore $\dist(A_1,A_2)<4\rho$, a contradiction.

The hypotheses $A_i\cap\gamma_i\neq\varnothing$ imply
$h_i\in\Lex{A_i}$.  Applying \eqref{eq:geometric-cotlar} with
$g=h_1$ and $h=h_2$, we obtain
\[
 q\in\Lthick{A_1}{\rho}
 \qquad\text{or}\qquad
 q^{-1}\in\Lthick{A_2}{\rho}.
\]
In the first case every geodesic from $o$ to $q.o$, in particular
$\gamma_3$, meets $B_\rho(a)$ for one fixed $a\in A_1$; in the second case
$\gamma_4$ meets $B_\rho(a)$ for one fixed $a\in A_2$.  This is exactly
\eqref{eq:set-long-product-alternative}.

\smallskip
\noindent
\textup{(i)}$\Rightarrow$\textup{(v)}.
Choose a hyperbolicity constant $\rho>0$ for $X$ and take $m_k=k$.  Fix the
data in \textup{(v)} and put $q=h_1h_2^{-1}$.  If
$|h_1.o|\geq|h_2.o|$, then
\[
 \frac{|q.o|+|h_1.o|-|h_2.o|}{2}
 \geq m_k+\rho
 \geq |a_1|+\rho.
\]
Since the geodesic $\gamma_1$ passes through $a_1$, Lemma
\ref{lem:branch-propagation}, applied to $(h_1,h_2^{-1})$, shows that every
geodesic from $o$ to $q.o$ meets $B_\rho(a_1)$.  In particular,
$\gamma_3\cap B_\rho(a_1)\neq\varnothing$.  If
$|h_2.o|>|h_1.o|$, the same argument applied to $(h_2,h_1^{-1})$ shows that
$\gamma_4\cap B_\rho(a_2)\neq\varnothing$.  Thus
\eqref{eq:point-long-product-alternative} holds.

\smallskip
\noindent
\textup{(v)}$\Rightarrow$\textup{(iv)}.
Choose
\[
 a_1\in A_1\cap\gamma_1,
 \qquad
 a_2\in A_2\cap\gamma_2.
\]
The pointwise alternative \eqref{eq:point-long-product-alternative} implies
\eqref{eq:set-long-product-alternative}, since
$B_\rho(a_i)\subset A_i^\rho$.

\smallskip
\noindent
\textup{(iv)}$\Rightarrow$\textup{(i)}.
By Lemma~\ref{lem:fixed-scaling}, it suffices to prove that, for every
nonprincipal ultrafilter $\omega$ and every basepoint sequence $(x_n)$, the
asymptotic cone
\[
 Y=\Cone_\omega(X,(x_n),(n))
\]
is an $\mathbb R$-tree.  Fix such a cone.  Suppose that $Y$ contains a
non-tripod internal geodesic triangle.  After trimming common initial
subsegments, we obtain a nontrivial simple internal triangle with vertices
$x,y,z$ and sides
\[
 \alpha=[x,y],\qquad \beta=[x,z],\qquad \gamma=[z,y].
\]
Choose a vertex sequence $(v_n)$ representing $x$.  By vertex transitivity,
choose $s_n\in\Gamma$ such that $s_n.v_n=o$.  Applying $s_n$
coordinatewise induces an isometry from $Y$ onto an asymptotic cone based at
the constant sequence $o$, sending $x$ to $o$.  We may therefore assume
$x=o$ and represent the other two vertices and the three internal sides by
geodesic triangles
\[
 \alpha_n=[o,g_n.o],\qquad
 \beta_n=[o,q_n.o],\qquad
 \gamma_n=[q_n.o,g_n.o]
\]
for suitable $g_n,q_n\in\Gamma$.

Parametrize $\alpha$ by arclength from $o$ and $\gamma$ by arclength from
$z$.  Choose
\[
 0<t<\min\left\{d(o,y),d(z,y),\frac12d(o,z)\right\}.
\]
Because the triangle is simple and $t/C_0>0$, the compact subsegments
$\alpha([t/C_0,t])$ and $\gamma([t/C_0,t])$ are disjoint from $\beta$.
Consequently,
\begin{equation}\label{eq:cone-positive-separation}
 \sigma:=\min\Big\{
   \dist\big(\alpha([t/C_0,t]),\beta\big),
   \dist\big(\gamma([t/C_0,t]),\beta\big)
 \Big\}>0.
\end{equation}

For $\omega$-almost every sufficiently large $n$, choose $k(n)$ so that
\[
 m_{k(n)}\leq tn<m_{k(n)+1}.
\]
The bounded-ratio condition gives
\begin{equation}\label{eq:scale-selection}
 \frac{t}{C_0}<\frac{m_{k(n)}}{n}\leq t.
\end{equation}
Put
\[
 h_n=q_n^{-1}g_n,
 \qquad
 \ell_n=m_{k(n)}.
\]
Because $\ell_n$ is an integer and the endpoints of the graph geodesics are
vertices, choose vertices $a_n\in\alpha_n$ and $c_n\in\gamma_n$ with
\[
 \dist(o,a_n)=\ell_n,
 \qquad
 \dist(q_n.o,c_n)=\ell_n,
\]
and set $b_n=q_n^{-1}.c_n$.  Let
\[
 \eta_n=q_n^{-1}.\gamma_n.
\]
Then $\eta_n$ is a geodesic from $o$ to $h_n.o$ passing through $b_n$, and
\begin{equation}\label{eq:cone-branches}
 |a_n|=|b_n|=\ell_n,
 \qquad
 g_nh_n^{-1}=q_n.
\end{equation}
Let
\[
 s=\lim_{n\to\omega}\frac{\ell_n}{n}\in[t/C_0,t].
\]
Then $(a_n)$ and $(c_n)$ represent the points $\alpha(s)$ and $\gamma(s)$,
respectively.  Hence \eqref{eq:cone-positive-separation} implies that, for
$\omega$-almost every $n$,
\begin{equation}\label{eq:cone-avoid}
 \dist(a_n,\beta_n)>10\rho,
 \qquad
 \dist(c_n,\beta_n)>10\rho.
\end{equation}
Moreover, $d(o,z)>2t$ and $\ell_n/n\leq t$, so, for $\omega$-almost every
$n$,
\begin{equation}\label{eq:cone-long}
 |q_n.o|\geq2\ell_n+2\rho.
\end{equation}
Let $\beta'_n$ be the reverse of the translated geodesic
$q_n^{-1}.\beta_n$.  Then $\beta'_n$ is a geodesic from $o$ to
$q_n^{-1}.o$, and
\begin{equation}\label{eq:cone-b'}
 \dist(b_n,\beta'_n)=\dist(c_n,\beta_n)>10\rho.
\end{equation}

Fix an index $n$ for which the preceding conclusions hold.  Apply
\textup{(iv)} at the permitted scale $\ell_n=m_{k(n)}$, with
\[
 h_1=g_n,\qquad h_2=h_n,\qquad
 \gamma_1=\alpha_n,\quad \gamma_2=\eta_n,\quad
 \gamma_3=\beta_n,\quad \gamma_4=\beta'_n.
\]
If $\dist(a_n,b_n)\geq4\rho$, take
$A_1=\{a_n\}$ and $A_2=\{b_n\}$.  The conclusion of \textup{(iv)} would
force
\[
 \beta_n\cap B_\rho(a_n)\neq\varnothing
 \qquad\text{or}\qquad
 \beta'_n\cap B_\rho(b_n)\neq\varnothing,
\]
contradicting \eqref{eq:cone-avoid} and \eqref{eq:cone-b'}.

If $\dist(a_n,b_n)<4\rho$, take
$A_1=A_2=A=\{a_n,b_n\}$.  By the triangle inequality,
\eqref{eq:cone-avoid} and \eqref{eq:cone-b'} imply that both $\beta_n$ and
$\beta'_n$ stay more than $6\rho$ from both points of $A$.  This contradicts
\eqref{eq:set-long-product-alternative} once again.

Thus every internal geodesic triangle in $Y$ is a tripod.  We next show that
geodesics in $Y$ are unique.  Fix an internal geodesic $[u,v]$, let $\tau$
be any other geodesic from $u$ to $v$, and choose $w\in\tau$.  Choose
internal geodesics from $u$ to $w$ and from $w$ to $v$.  Together with
$[u,v]$, they form an internal tripod.  Since
\[
 d(u,v)=d(u,w)+d(w,v),
\]
its branch point must be $w$, and hence $w\in[u,v]$.  Thus
$\tau=[u,v]$, so $Y$ is an $\mathbb R$-tree.  Since $\omega$ and $(x_n)$
were arbitrary, Lemma~\ref{lem:fixed-scaling} implies that $X$ is
hyperbolic.
\end{proof}

\begin{corollary}\label{cor:cayley-cotlar-characterization}
A finitely generated group $\Gamma$ is hyperbolic if and only if one
(equivalently, every) Cayley action has the $\rho$-Cotlar property for some
$\rho>0$.  In the converse direction, it is enough that the
$(\rho,m)$-Cotlar property hold along an increasing sequence of radii
satisfying \eqref{eq:bounded-ratio-scales}.
\end{corollary}

\begin{remark}
The conclusion of Theorem~\ref{thm:cotlar-characterization} concerns the
geometry of the vertex-transitive graph $X$, not that of $\Gamma$.  Freeness,
faithfulness, and properness of the action are not required.  For example,
$\mathbb Z^2$ acts transitively on the line $\mathbb Z$ through the first
coordinate; the graph is a tree even though $\mathbb Z^2$ is not a hyperbolic
group.

The bounded-ratio condition may equivalently be written as follows: there is
$C\geq1$ such that every sufficiently large interval $[R/C,R]$ contains a
permitted radius.  Thus the converse only requires the Cotlar property on a
set of scales having bounded gaps on the logarithmic scale.
\end{remark}

\section{A bounded Gromov interface criterion}
\label{sec:bounded-interface}

The preceding results use a thin/thick pair of geometric branch symbols in a hyperbolic space. We now record a complementary exact, unthickened Cotlar criterion for more general partitions. No hyperbolicity of the ambient metric space is assumed in this section.
Throughout this section, $X$ is an arbitrary metric space and
$\Gamma\curvearrowright X$ is an action by isometries.  We put
\[
 \ell_X(g)=\dist(o,g.o),\qquad g\in\Gamma.
\]

Given a pairwise disjoint family $(\mathcal E_i)_{i\in I}$ of subsets of
$\Gamma$, after relabeling if necessary assume that $0\notin I$, and put
\[
 \mathcal E_0=\Gamma\setminus\bigcup_{i\in I}\mathcal E_i,
 \qquad I_0=I\cup\{0\}.
\]
Thus $(\mathcal E_i)_{i\in I_0}$ is a partition of $\Gamma$.  We formulate
the criterion directly for such partitions.

\begin{definition}[Bounded Gromov interface with respect to the action]
\label{def:bounded-interface}
Let $(\mathcal E_i)_{i\in I_0}$ be a partition of $\Gamma$.  We say that
this partition has a \emph{bounded Gromov interface with respect to the
action $\Gamma\curvearrowright X$} if there exists $M\geq0$ such that
\begin{equation}\label{eq:bounded-interface}
 g\in\mathcal E_i,\quad h\in\mathcal E_j,\quad i\neq j
 \quad\Longrightarrow\quad
 (g.o\mid h.o)_o\leq M.
\end{equation}
\end{definition}

Thus orbit points belonging to different parts of the partition cannot
fellow travel from $o$ for an arbitrarily long time.  The constant $M$ is
uniform over the whole partition.

\begin{theorem}[Bounded-interface multiplier theorem]
\label{thm:bounded-interface}
Let $(\mathcal E_i)_{i\in I_0}$ be a partition of $\Gamma$ with bounded
Gromov interface constant $M$.  Choose $R>2M$ and assume that $E_{<R}$ is
finite.  For every   family
$\eps=(\eps_i)_{i\in I_0}$ with $|\eps_i|\leq1$, define
\[
 T_\eps=\sum_{i\in I_0}\eps_iP_{\mathcal E_i},
\]
as a map on $L^2(\hat\Gamma)$. Then, for every $1<p<\infty$, there is a constant $C_p$ such that
\begin{equation}\label{eq:bounded-interface-Lp}
 \|T_\eps x\|_p\leq C_p\|x\|_p,
 \qquad x\in L^p(\widehat\Gamma).
\end{equation}
  Consequently, the subspaces
\[
 X_i^p=\overline{\operatorname{span}}^{\,L^p}
       \{\lambda_g:g\in\mathcal E_i\},
 \qquad i\in I_0,
\]
form an $L^p$-unconditional decomposition of
$L^p(\widehat\Gamma)$.
\end{theorem}

\begin{proof}
For each $g\in\Gamma$, let $i(g)\in I_0$ be the unique index satisfying
$g\in\mathcal E_{i(g)}$, and put
\[
 m_\eps(g)=\eps_{i(g)},
 \qquad
 T_\eps(\lambda_g)=m_\eps(g)\lambda_g.
\]
Since $|m_\eps(g)|\leq1$, the multiplier $T_\eps$ is an
$L^2$-contraction.  We first derive a scalar Cotlar identity.  Let
$g,h\in\Gamma$ and put $q=gh^{-1}$.  If both
\[
 m_\eps(g)-m_\eps(q)
 \quad\text{and}\quad
 m_\eps(h)-m_\eps(q^{-1})
\]
are nonzero, then $g$ and $q$ belong to different parts of the partition,
and so do $h$ and $q^{-1}$.  Hence \eqref{eq:bounded-interface} gives
\[
 \frac{\ell_X(g)+\ell_X(q)-\ell_X(h)}2
   =(g.o\mid q.o)_o\leq M,
\]
because $\dist(g.o,q.o)=\ell_X(h)$, and similarly
\[
 \frac{\ell_X(h)+\ell_X(q)-\ell_X(g)}2
   =(h.o\mid q^{-1}.o)_o\leq M.
\]
Adding these inequalities yields $\ell_X(q)\leq2M$.  Therefore, whenever
$q\in E_{\geq R}$,
\begin{equation}\label{eq:interface-factor}
 \big(m_\eps(g)-m_\eps(q)\big)
 \overline{\big(m_\eps(h)-m_\eps(q^{-1})\big)}=0.
\end{equation}
Expanding \eqref{eq:interface-factor}, we obtain
\begin{equation}\label{eq:interface-scalar-cotlar}
\begin{aligned}
 m_\eps(g)\overline{m_\eps(h)}
 ={}&m_\eps(q)\overline{m_\eps(h)}
   +m_\eps(g)\overline{m_\eps(q^{-1})}\\
  &-m_\eps(q)\overline{m_\eps(q^{-1})}.
\end{aligned}
\end{equation}

Define the right-handed multiplier by
\[
 T_\eps^\circ(\lambda_s)
   =\overline{m_\eps(s^{-1})}\lambda_s.
\]
Then $(T_\eps x)^*=T_\eps^\circ(x^*)$, and
\eqref{eq:interface-scalar-cotlar} gives, for every Fourier polynomial $x$,
\begin{equation}\label{eq:interface-operator-cotlar}
\begin{aligned}
 P_{E_{\geq R}}\big(T_\eps x\,T_\eps^\circ(x^*)\big)
 =P_{E_{\geq R}}\Big(&
 T_\eps\big(xT_\eps^\circ(x^*)\big)
 +T_\eps^\circ\big(T_\eps(x)x^*\big)\\
 &-T_\eps^\circ T_\eps(xx^*)\Big).
\end{aligned}
\end{equation}

Taking $L^2$-norms in \eqref{eq:interface-operator-cotlar} and repeating
the low--high frequency argument from the proof of
Theorem~\ref{thm:main} gives the $L^4$ estimate, uniformly over
$\|\eps\|_\infty\leq1$.  The same dyadic induction, followed by
interpolation and duality, proves \eqref{eq:bounded-interface-Lp} for
every $1<p<\infty$.
 \end{proof}

\begin{corollary}[Compact commutators for bounded-interface partitions]
\label{cor:interface-commutators}
Let $(\mathcal E_i)_{i\in I_0}$ be a partition of $\Gamma$ with bounded
Gromov interface constant $M$, and assume that $E_{<R}$ is finite for every
$R>0$.  For $i\in I_0$, put
\[
 P_i^\circ=P_{\mathcal E_i^{-1}}.
\]
Then $P_i^\circ$ is bounded on $L^p(\widehat\Gamma)$ for every
$1<p<\infty$.  Moreover, for every $r>2$ and every
$y\in L^r(\widehat\Gamma)$, the commutator
\[
 [P_i^\circ,L_y]:\mathcal L(\Gamma)\longrightarrow L^2(\widehat\Gamma)
\]
is compact.
\end{corollary}

\begin{proof}
Taking the coefficient family supported only at $i$ in
Theorem~\ref{thm:bounded-interface} gives the $L^p$-boundedness of
$P_{\mathcal E_i}$.  Since
\[
 P_i^\circ(x)=\big(P_{\mathcal E_i}(x^*)\big)^*,
\]
the same conclusion holds for $P_i^\circ$.

Fix $s,h\in\Gamma$.  On the Fourier basis,
\[
 \big(P_i^\circ L_{\lambda_s}-L_{\lambda_s}P_i^\circ\big)
 (\lambda_{h^{-1}})
 =\big(\mathbf 1_{\mathcal E_i}(hs^{-1})
       -\mathbf 1_{\mathcal E_i}(h)\big)\lambda_{sh^{-1}}.
\]
If the coefficient on the right is nonzero, then $h$ and $hs^{-1}$ lie in
different parts of the partition.  Hence
\[
 (h.o\mid hs^{-1}.o)_o\leq M.
\]
On the other hand,
\begin{align*}
 (h.o\mid hs^{-1}.o)_o
 &=\frac{\ell_X(h)+\ell_X(hs^{-1})-\ell_X(s)}2\\
 &\geq \ell_X(h)-\ell_X(s).
\end{align*}
Thus $\ell_X(h)\leq M+\ell_X(s)$.  Since every orbit ball is finite, only
finitely many $h$ can occur, and
$[P_i^\circ,L_{\lambda_s}]$ has finite rank.  The same is true when $y$ is
a Fourier polynomial.

Let $r>2$ and choose $q$ so that
\[
 \frac12=\frac1r+\frac1q.
\]
For $x\in\mathcal L(\Gamma)$, H\"older's inequality and the
$L^q$-boundedness of $P_i^\circ$ give
\begin{align*}
 \|P_i^\circ(yx)-yP_i^\circ(x)\|_2
 &\leq \|yx\|_2+\|yP_i^\circ(x)\|_2\\
 &\leq \big(1+\|P_i^\circ\|_{L^q\to L^q}\big)
          \|y\|_r\|x\|_\infty.
\end{align*}
Therefore the commutator depends continuously on $y$ in the $L^r$-norm.
Approximating $y$ by Fourier polynomials proves compactness.
\end{proof}

\subsection{Examples and sources of bounded Gromov interfaces}
\label{subsec:interface-examples}

The bounded-interface condition is naturally produced by continuous
finite-valued labels on a suitable compactification of the orbit.  We say
that a compact metrizable compactification
\[
 \overline O=O\sqcup\partial O,
 \qquad O=\Gamma.o,
\]
has the \emph{Gromov convergence property} if, whenever
$x_n\to\xi$, $y_n\to\eta$ in $\overline O$ and
$(x_n\mid y_n)_o\to\infty$, one has
$\xi=\eta\in\partial O$.  For a proper geodesic hyperbolic space, the
usual Gromov compactification has this property.  In the nonproper setting,
we only assume the stated compactification of the orbit.

\begin{example}[Continuous finite-valued labels]
\label{ex:continuous-labels}
Assume that $X$ is a geodesic hyperbolic metric space and that the orbit
$O=\Gamma.o$ admits a compact metrizable compactification
$\overline O=O\sqcup\partial O$ with the Gromov convergence property.  Let
\[
 \Gamma=\bigsqcup_{i=0}^N\mathcal E_i
\]
be a finite partition saturated with respect to the orbit map.  Equivalently,
the label map
\[
 c:O\longrightarrow\{0,1,\ldots,N\},
 \qquad c(g.o)=i\quad\text{when }g\in\mathcal E_i,
\]
is well defined.  Suppose that $c$ extends continuously to $\overline O$,
where the target is given the discrete topology.  Then the partition has a
bounded Gromov interface.

Indeed, otherwise, after passing to subsequences, there would be fixed
$i\neq j$ and sequences $g_n\in\mathcal E_i$, $h_n\in\mathcal E_j$ such
that
\[
 (g_n.o\mid h_n.o)_o\longrightarrow\infty.
\]
Compactness and the Gromov convergence property would force both orbit
sequences to converge to the same point of $\partial O$, contradicting the
continuity of $c$.  Consequently, whenever the low-frequency hypothesis of
Theorem~\ref{thm:bounded-interface} is satisfied, the associated Fourier
subspaces form an $L^p$-unconditional decomposition of
$L^p(\widehat\Gamma)$ for every $1<p<\infty$.  If every orbit ball is
finite, Corollary~\ref{cor:interface-commutators} also gives compact
commutators for the right Fourier projection associated with every part.
\end{example}

\begin{example}[Finite-range continuous boundary functions]
\label{ex:finite-range-boundary-symbols}
Let $\Gamma$ be a finitely generated hyperbolic group,
let $X=\operatorname{Cay}(\Gamma,S)$, and write
\[
 \overline\Gamma=\Gamma\sqcup\partial\Gamma
\]
for the Gromov compactification.  Suppose that
$f\in C(\partial\Gamma)$ has finite range
\[
 f(\partial\Gamma)=\{c_0,\ldots,c_N\}.
\]
Then $f$ admits a finite-range continuous extension
$\widetilde f\in C(\overline\Gamma)$.  For every such extension, the
restriction
\[
 m=\widetilde f|_\Gamma
\]
defines a bounded Fourier multiplier
\[
 T_m:L^p(\widehat\Gamma)\longrightarrow L^p(\widehat\Gamma),
 \qquad 1<p<\infty.
\]

To construct an extension, put
\[
 F_i=f^{-1}(\{c_i\}),\qquad 0\leq i\leq N.
\]
The sets $F_i$ form a finite clopen partition of $\partial\Gamma$.
Since $\overline\Gamma$ is compact Hausdorff, they admit pairwise
disjoint open neighborhoods $U_i\subset\overline\Gamma$.  Their union
contains $\partial\Gamma$, so
\[
 K=\overline\Gamma\setminus\bigcup_{i=0}^N U_i
\]
is a compact subset of the discrete open set $\Gamma$ and is therefore
finite.  Define $\widetilde f=c_i$ on
$F_i\cup(U_i\cap\Gamma)$ and assign to the finitely many points of $K$
arbitrary values from $\{c_0,\ldots,c_N\}$.  This gives the desired
finite-range continuous extension.

For $0\leq i\leq N$, let
\[
 \mathcal E_i=\{g\in\Gamma:m(g)=c_i\}.
\]
The level-set indicator $\mathbf 1_{\mathcal E_i}$ extends continuously
to $\overline\Gamma$; for example, it is the restriction of
\[
 \prod_{j\neq i}
 \frac{\widetilde f-c_j}{c_i-c_j}.
\]
Hence $(\mathcal E_i)_{i=0}^N$ is a boundary-continuous finite partition,
so Example~\ref{ex:continuous-labels} gives a bounded Gromov interface.
Since word balls are finite, Theorem~\ref{thm:bounded-interface} applies,
and
\[
 T_m=\sum_{i=0}^N c_iP_{\mathcal E_i}
\]
is bounded on every $L^p(\widehat\Gamma)$, $1<p<\infty$.

This multiplier is canonical modulo finite Fourier support.  Indeed, if
$\widetilde f_1$ and $\widetilde f_2$ are two finite-range continuous
extensions of the same boundary function, then
$\widetilde f_1-\widetilde f_2$ vanishes on $\partial\Gamma$.  Its
restriction to $\Gamma$ belongs to $c_0(\Gamma)$ and has finite range,
so it has finite support.  The corresponding multipliers therefore differ
by a finite-rank Fourier multiplier.

The finite-range hypothesis is essential for this argument.  Since
$\partial\Gamma$ is compact, every continuous boundary function is bounded,
but the present example does not assert that an arbitrary continuous
scalar boundary function gives an $L^p$-multiplier.
\end{example}

\begin{example}[Exact branches and free-group prefix partitions]
\label{ex:prefix-interface}
For exact branch sets $\mathcal E_i=\Lex{A_i}$, saturation with respect to
the orbit map is automatic.  Thus, if finitely many pairwise disjoint exact
branches and their complement define a label on $O$ which extends
continuously to an orbit compactification as in
Example~\ref{ex:continuous-labels}, then they form a bounded-interface
partition.  Notice that the required continuity concerns the exact branch
label on the whole orbit compactification, not merely the clopenness of a
subset of $\partial O$ without compatibility with its values on $O$.

The standard Cayley tree of a finitely generated free group gives the basic
model.  Let $\Gamma=\mathbb F_n$, fix $m\geq1$, and let
\[
 S_m(e)=\bigsqcup_{i\in I}A_i
\]
be any partition.  Put
\[
 \mathcal E_0=E_{<m},
 \qquad
 \mathcal E_i=\Lex{A_i}\quad(i\in I).
\]
Then $(\mathcal E_i)_{i\in I\cup\{0\}}$ is a partition of $\mathbb F_n$
with bounded Gromov interface constant $m-1$.  Indeed, in a rooted tree the
Gromov product is the length of the common initial segment.  If two words
belong to different length-$m$ prefix branches, their common initial segment
has length at most $m-1$; if one belongs to $E_{<m}$, the same bound is
immediate.

In particular, for $a\in S_m(e)$ one may take
\[
 \mathcal E_1=\Lex{\{a\}},
 \qquad
 \mathcal E_2=\Lex{S_m(e)\setminus\{a\}},
 \qquad
 \mathcal E_0=E_{<m}.
\]
Theorem~\ref{thm:bounded-interface} gives the corresponding global
$L^p$-unconditional decomposition, and
Corollary~\ref{cor:interface-commutators} gives compact commutators for the
right branch projections.  For $m=1$ and the singleton partition of
$S_1(e)$, these statements recover the free-group results of Mei--Ricard
\cite{MR17}, including the positive answer to Ozawa's question from
\cite{O10}.
\end{example}

\begin{example}[A hyperbolic-group non-example]
\label{ex:unbounded-exact-interface}
Hyperbolicity alone does not force the natural exact
existential/universal partition to have bounded interface.  Let
\[
 \Gamma=\mathbb F(u,v)\times\mathbb Z_2
\]
and let $c$ denote the nontrivial element of the central $\mathbb Z_2$
factor.  Equip $\Gamma$ with the generating set
$\{u^{\pm1},v^{\pm1},c\}$.  This is a non-elementary hyperbolic group,
since it is a finite extension of $\mathbb F(u,v)$.  Fix $m\geq1$, put
\[
 a=u^{m-1}c\in S_m(e),
\]
and consider
\[
 \mathcal E_0=E_{<m},
 \qquad
 \mathcal E_1=\Lex{\{a\}},
 \qquad
 \mathcal E_2=\Lall{S_m(e)\setminus\{a\}}.
\]
These three sets form a partition of $\Gamma$.  Indeed, an element of length
at least $m$ either has a geodesic passing through $a$, or every one of its
geodesics meets $S_m(e)\setminus\{a\}$.

For $N\geq m$, set
\[
 g_N=u^Nc,
 \qquad
 h_N=u^N.
\]
The geodesic word
$u^{m-1}cu^{N-m+1}$ shows that $g_N\in\mathcal E_1$.  Moreover, the
unique geodesic from $e$ to $h_N=u^N$ is labelled by $u^N$; it meets
$S_m(e)$ at $u^m\neq a$.  Hence $h_N\in\mathcal E_2$.  On the other hand,
\[
 |g_N|=N+1,
 \qquad |h_N|=N,
 \qquad \dist(g_N,h_N)=1,
\]
and hence
\[
 (g_N\mid h_N)_e=N\longrightarrow\infty.
\]
Thus this exact existential/universal partition does not have a bounded
Gromov interface.  Equivalently, the exact branch label fails the
asymptotic stability expressed by Example~\ref{ex:continuous-labels}.
\end{example}

\section*{Acknowledgements}
The author thanks S. Kunnawalkam Elayavalli, M. Kalantar, N. Ozawa,
E. Ricard, and R. Toyota for helpful discussions.  The author also thanks
M. Junge, H. Lee, and J. Parcet for encouraging the development of a
version of the Mei--Ricard theory for groups acting on hyperbolic spaces.

\end{document}